\documentclass[11pt]{amsart}

\usepackage{palatino}
\usepackage[all]{xy,xypic}
\usepackage{amsfonts,amsmath,oldgerm,amssymb,amscd,tikz-cd,enumerate}
\usepackage[a4paper,left=12em,right=12em,top=10em,bottom=10em]{geometry}

\newcommand{\ra}{\rightarrow}		

\newcommand{\by}[1]{\stackrel{#1}{\ra}}

\newcommand{\remove}[1]{}

\newcommand{\inj}{\hookrightarrow}
\newcommand{\ol}{\overline}		

\newcommand{\iso}{\by \sim}

\newtheorem{theorem}{Theorem}[section]

\newtheorem{proposition}[theorem]{Proposition}
\newtheorem{lemma}[theorem]{Lemma}
\newtheorem{definition}[theorem]{Definition}
\newtheorem{corollary}[theorem]{Corollary}

\newtheorem{conjecture}[theorem]{Conjecture}
\newtheorem{example}[theorem]{Example}

	\newcommand{\CP}{\mbox{$\mathcal P$}}
	
\newcommand{\CS}{\mbox{$\mathcal S$}}

\newcommand{\MM}{\mbox{$\mathfrak M$}}	\newcommand{\MN}{\mbox{$\mathfrak N$}}

\newcommand{\MS}{\mbox{$\mathfrak S$}}	\newcommand{\MT}{\mbox{$\mathfrak T$}}

\newcommand{\mm}{\mbox{$\mathfrak m$}}	
	\newcommand{\p}{\mbox{$\mathfrak p$}}
\newcommand{\mq}{\mbox{$\mathfrak q$}}

\newcommand{\Spec}{\text{Spec}}	
\newcommand{\hh}{\text{ht}}
\newcommand{\rank}{\text{rank}}

\newcommand{\Aut}{\mbox{\rm Aut\,}}

\newcommand{\Um}{\text{Um}}		
\newcommand{\SL}{\text{SL}}
\newcommand{\E}{\text{E}}	
\newcommand{\GL}{\text{GL}}

\newcommand{\bp}{\begin{proposition}}
	\newcommand{\ep}{\end{proposition}}
\newcommand{\bl}{\begin{lemma}}
	\newcommand{\el}{\end{lemma}}
\newcommand{\bt}{\begin{theorem}}
	\newcommand{\et}{\end{theorem}}
\newcommand{\bc}{\begin{corollary}}
	\newcommand{\ec}{\end{corollary}}
\newcommand{\bd}{\begin{definition}}
	\newcommand{\ed}{\end{definition}}
\newcommand{\bco}{\begin{conjecture}}
	\newcommand{\eco}{\end{conjecture}}
\newcommand{\bma}{\begin{bmatrix}}
	\newcommand{\ema}{\end{bmatrix}}

\def\rmk{\refstepcounter{theorem}\paragraph{{\bf Remark} \thetheorem}}
\def\proof{\paragraph{Proof}}
\def\example{\refstepcounter{theorem}\paragraph{{\bf Example} \thetheorem}}
\def\notation{\paragraph{\bf Notation}}
\def\quest{\refstepcounter{theorem}\paragraph{{\bf Question} \thetheorem}}

\title [Stability theorems for positively graded domains and a question of Lindel]{Stability theorems for positively graded domains\\ and a question of Lindel}
\author{Sourjya Banerjee}
\address{The Institute of Mathematical Sciences, HBNI, C.I.T. Campus, Tharamani, Chennai  600113, India}
\email{sourjyab@imsc.res.in, sourjya91@gmail.com}
\date{\today}


\usepackage{hyperref}

\begin{document}
	\maketitle

\begin{abstract}
Given a commutative Noetherian graded domain $R = \bigoplus_{i\ge 0} R_i$ of dimension $d\geq 2$ with $\dim(R_0) \geq 1$, we prove that any unimodular row of length $d+1$ in $R$ can be completed to the first row of an invertible matrix $\alpha$ such that $\alpha$ is homotopic to the identity matrix. Utilizing this result, it has been established that if $I \subset R$ is an ideal satisfying $\mu(I/I^2) = \hh(I) = d$, then any set of generators of $I/I^2$ lifts to a set of generators of $I$, where $\mu(-)$ denotes the minimal number of generators. Consequently, any projective $R$-module of rank $d$ with trivial determinant splits into a free factor of rank one. This provides an affirmative answer to an old question of Lindel. Finally, we prove that for any projective $R$-module $P$ of rank $d$, if the Quillen ideal of $P$ is non-zero, then $P$ is cancellative.
	
\end{abstract}

\section{Introduction}
We commence by recalling an old question of Murthy \cite{Bass73}. Let $A=\bigoplus_{i\ge 0} R_i$ be a normal positively graded finitely generated algebra over $R_0$, where $R_0=k$ is a field. Then Murthy asked whether $\text{K}_0(A) \cong \mathbb{Z}$. Bloch provided a counterexample to this question by considering $A=\frac{\mathbb{C}[X,Y,Z]}{\langle Z^2 - X^3 - Y^7\rangle}$. However, if (1) $k$ is an algebraically closed field of characteristic $p> 0$, (2) $A$ is a Cohen-Macaulay ring of dimension $2$, and (3) the vertex (corresponding to the ideal $R_{+}:=\bigoplus_{i\ge 1}R_i$) is the only singularity of $\Spec (A)$, then Srinivas \cite[Corollary 1.3]{sv1} showed that Murthy's question has an affirmative answer. Therefore, using the cancellation theorem of Murthy and Swan \cite{PMRS}, it follows that every projective $A$-module is free. Hence, this improves the existing stability theorems for projective modules over such graded algebras of dimension $2$. 

In 1987, Lindel \cite[Theorem 1.3]{Li87} improved Quillen's Local-Global Principle from polynomial rings to positively graded rings. This, in particular, initiated the study of projective modules over a higher dimensional graded ring from a more algebraic point of view. Let us briefly recall Lindel's work from \cite{Li87}. Consider $A=\bigoplus_{i\ge 0} R_i$, a commutative Noetherian (non-trivially) graded ring of dimension $d\ge 2$, and let $M$ be a finitely presented module over $A$. Lindel first established in \cite[Theorem 1.3]{Li87} that the Quillen set of $M$, defined by $J(R_0,M):=\{s\in R_0:M_s\cong \frac{M_s}{(MR_+)_s}\otimes R_s\}$, forms an ideal in $R_0$. Furthermore, in the same article, he proved the following theorem.

\bt\cite[Theorem 2.5]{Li87}\label{lt}
Let $A=\bigoplus_{i\ge 0} R_i$ be a commutative Noetherian (non-trivially) graded ring of dimension $d\ge 2$ and let $P$ be a projective $A$-module of rank $d$ such that $A$ and $P$ satisfy the following conditions.

\begin{enumerate}
	\item $\dim(R_0)=d-1$;
	\item $A=R_0[t_1,\cdots, t_n]$, where $t_i$ are homogeneous in $R_{+}$ for $i=1,\cdots, n$ such that the kernel of the $R_0$-epimorphism $\phi: R_0[T_1,\cdots,T_n]\twoheadrightarrow R$ sending $T_i\mapsto t_i$ has a height $\ge n-1$;
	\item $P_{1+JR_+}$ splits into a free summand of rank one, where $J=J(R_0,P)$.
\end{enumerate} 
Then $P$ splits into a free summand of rank one.
\et

In \cite[Remark 2.6]{Li87}, Lindel queries whether the hypothesis (2) in his theorem is necessary. In this article we prove that when $R$ is a positively graded domain, all the hypotheses (1) to (3) in Theorem \ref{lt} are redundant. Specifically, we establish the following much stronger version, the proof of which can be found in Theorems \ref{eg}, \ref{egd2}, and Corollary \ref{eue}.

\bt\label{egi}
Let $A=\bigoplus_{i\ge 0} R_i$ be a commutative Noetherian (non-trivially) graded domain of dimension $d\ge 2$ such that $\dim(R_0)\ge 1$. Let $C$ and $n$ be one of the following:
\begin{enumerate}
	\item $C=A$ and $n=\dim(A)=d\ge 2$.
	\item $C=S^{-1}A$, where $S\subset A$ is a multiplicative set contained in the set of all non-zero divisors in $A$ such that $\dim(C)=\dim(A)$ and $n=\dim(A)=d\ge 3$.
	\item $C=B[M]$, where $B$ is a commutative Noetherian ring of dimension $ \ge 2$ and $M$ is a finitely generated commutative cancellative (not necessarily torsion free) monoid of rank $r\ge 1$. We take $n=\dim(B[M])$.
\end{enumerate}
Let $I\subset C$ be an ideal such that $\mu(I/I^2)=\hh(I)=n$, where $\mu(-)$ denotes the minimal number of generators. Then any set of generators of $I=\langle f_1,\ldots,f_n\rangle +I^2$ lifts to a set of generators of $I$. Consequently, any projective $C$-module of rank $n$ (with trivial determinant) splits into a free summand of rank one. 
\et

Readers may question the significance of the hypothesis $\dim(R_0)\ge 1$ in this article. However, in Examples \ref{ex1} and \ref{ex2}, we demonstrate that this hypothesis is indeed necessary in Theorem \ref{egi}. On the other hand, to improve \cite[Theorem 2.5]{Li87}, we needed to take a significantly distinct approach from Lindel's. A crucial step in establishing Theorem \ref{egi} are Theorem \ref{umch} and Corollary \ref{ec1}, where we prove the following result.

\bt\label{hc}
Let $A$ be as in Theorem \ref{egi}. Let $C=S^{-1}A$, where $S\subset A$ is a multiplicative set. Any unimodular row in $C$ of length $d+1$ can be completed to the first row of an invertible matrix $\alpha$. Moreover, the matrix $\alpha$ can be chosen in such a way that there exists $\theta(T)\in \GL_{d+1}(C[T])$ such that $\theta(0)=\text{Id}$ and $\theta(1)=\alpha$ (in this case $\alpha$ will be called homotopic to the identity matrix). Consequently, any stably free $C$-module of rank $d$ is free.
\et

It may be observed that Theorem \ref{hc}, over an arbitrary ring, does not imply Theorem \ref{egi} (1) and (2). An example illustrating this point is provided in \cite[Remark 3.8]{NMR}. To establish such implications within our framework, we utilize the additional graded structure of the ring. Even to prove Theorem \ref{hc}, the techniques used in this article [e.g., Proposition \ref{splitting} and Lemma \ref{patching}] are very specific to the graded set-up, and do not extend to arbitrary rings [see Remark \ref{rpvul}]. Furthermore, we generalize Theorem \ref{hc} in the following form [for the proof, we refer to Theorem \ref{can1}].

\bt
Let $A$ be as in Theorem \ref{egi}. Let $P$ be a projective $R$-module of rank $d$ such that $J(R_0,P)\not=0$. If $P\oplus A^k\cong Q\oplus A^k$ for some $k\in \mathbb{N}$, then $P\cong Q$. In other words, the module $P$ is cancellative. 
\et

\subsection{On a question of Nori over a graded non-smooth algebra} In Section \ref{a1}, we deduce some consequences of Theorem \ref{egi}. We very briefly recall an algebraic analogy of a question asked by Nori \cite{SM}. 
\smallskip

\quest\label{Norifree}
Let $C$ be a smooth affine domain of dimension $d$ over an infinite perfect field. Let $I\subset C[T]$ be an ideal of height $n$ such that $I=\langle f_1,\ldots,f_n\rangle +I^2T$, where $2n\geq d+3$. Do there exist $g_i\in I$ such that $I=\langle g_1,\ldots,g_n\rangle $ with $g_i-f_i\in I^2T$?
   
This question is completely solved in \cite{BR} and \cite{BHMK}. Bhatwadekar, Mohan Kumar and Srinivas constructed an example \cite[Example 6.4]{BR} of a non-smooth positively graded affine domain (with the degree zero subring $ \mathbb{C}$) such that over which Nori's question has a negative answer. However, when the ring has singularities, it is shown in \cite{sbmkd1} that imposing some suitable smoothness condition on the ideal $I\cap R$ one can prevent such anomalies. Here, in Section \ref{a1}, we aim to understand the underlying issue that prevents the existence of such a lift in \cite[Example 6.4]{BR}. In particular, we prove the following [for details we refer to Theorem \ref{NQ} and Corollary \ref{egpa}].

\bt
Let $A=\bigoplus_{i\ge 0} R_i$ be an affine domain (non-necessarily smooth) of dimension $d\ge 3$ over an infinite field such that $\frac{1}{d!}\in A$ and $\dim(R_0)\ge 1$. Let $I\subset A[T]$ an ideal such that $\mu(I/I^2T)=\hh(I)=d$. Then any set of generators of $I/I^2T$ lifts to a set of generators of $I$. Consequently, any projective $A[T]$-module (with trivial determinant) of rank $d$ splits into a free summand of rank one.
\et
\subsection{{Layout of the article}}The article is organized as follows: Section \ref{2} covers basic definitions and preliminary results necessary for proving the remaining parts of the article. In Section \ref{ur}, we present various results related to the completion of unimodular rows. The main result of this section is Theorem \ref{umch}. In this section we also improve the existing injective stability bounds for the classical group $\mathrm{SK_1}$ of graded domains (see Theorem \ref{ist}). Section \ref{MT} is dedicated to proving the main theorems of this article, specifically Theorem \ref{eg}. Section \ref{cansection} addresses the cancellation property for projective modules over graded domains. Finally, we conclude with some applications in Section \ref{a1}.

\subsection{{Convention}} The symbol $\mathbb{N}$ denotes the set of all non-negative integers, including $0$. All rings considered in this article are assumed to be commutative Noetherian with $1(\neq 0)$ having finite (Krull) dimension. Additionally, all graded rings discussed in this article are assumed to have a non-trivial $\mathbb{N}$-grading. For a graded ring $R=\bigoplus_{i\geq 0} R_i$, we use the notation $R_+=\bigoplus_{i\geq 1} R_i$ to represent the irrelevant ideal in $R$ containing all elements which can be written as a sum of homogeneous elements of degree $>0$. Every module considered in this article is assumed to be finitely generated. The symbol $e_1$ denotes the vector $(1, 0, \ldots, 0)$.

\section{Preliminaries}\label{2}
\renewcommand{\bt}{\begin{theorem}}
	\renewcommand{\et}{\end{theorem}}
\renewcommand{\bc}{\begin{corollary}}
	\renewcommand{\ec}{\end{corollary}}
\def\rmk{\refstepcounter{theorem}\paragraph{{\bf Remark} \thetheorem}}
This section summarizes several results and definitions from the literature that are frequently used in this article to prove the main theorems. We may restate or improve these results as necessary. Before proceeding further, we recall several definitions from the literature.


\bd\label{def} Let $A$ be a ring.
\begin{enumerate}  
	\item Let $M$ be an $A$-module. An element $x\in M$ is said to be a \textit{basic element} of $M$ at a prime ideal $\p\in \Spec(A)$ if $x\not \in \p M_{\p}.$ For any $\CS\subset \Spec(A)$, we call $x$ a basic element of $M$ on $\CS$ if it is a basic element of $M$ at each prime ideal $\p\in \CS$.
	\item Let $\CS\subset \Spec(A)$ and let $\delta:\CS\to \mathbb{N}$ be a function. For two prime ideals $\p,\mq\in \CS$, we define a partial order $\p<<\mq$ if and only if $\p\subset \mq$ and $\delta(\p)>\delta(\mq)$.  We say that $\delta$ is a \textit{generalized dimension function} if for any ideal $I\subset A$, the set $V(I)\cap \CS$ has only finitely many minimal elements with respect to $<<$.
	\item A row vector $v = (v_1, \ldots, v_n) \in A^n$ is called a unimodular row of length $n$ if there exists $(\lambda_1, \ldots, \lambda_n) \in A^n$ such that $\lambda_1v_1 + \ldots + \lambda_n v_n = 1$. The set of all unimodular rows in $A^n$ of length $n$ is denoted by $\Um_n(A)$.
	\item  A positive integer $r$ is said to be the \textit{stable rank} of $A$, denoted as $\text{sr}(A)$, if $r$ is the smallest integer for which any $(a_1,\ldots,a_{r+1})\in \Um_{r+1}(A)$, there exist $\lambda_i\in A$, $i=1,\ldots,r$ such that $(a_1+\lambda_1a_{r+1},\ldots,a_r+\lambda_ra_{r+1})\in \Um_r(A)$.
	\item Let $\lambda \in A$, and let $e_{ij}(\lambda)$ denote the $n\times n$ matrix whose only possible non-zero entry is $\lambda$ at the position $(i,j)$, where $n\in \mathbb N$. We define $\E_n(A)$ as the subgroup of $\SL_n(A)$ generated by the matrices $E_{ij}(\lambda) := \text{Id} + e_{ij}(\lambda)$, where $\lambda \in A$ and $i \neq j$.
	\item Let $\alpha \in \E_n(A)$. Then $\alpha$ can be viewed as an element of $\E_{n+1}(A)$ via the canonical inclusion $\alpha \inj \begin{pmatrix}
		\alpha & 0 \\
		0 & 1
	\end{pmatrix}$. We define $\E(A) := \bigcup_{i \in \mathbb{N}} \E_i(A)$.
	\item We define $\text{H}_n(A)=\{\alpha\in \GL_n(A):\text{ there exists a } \theta(T)\in \GL_{n}(A[T]) \text{ such that } \theta(0)=\text{Id}$ \text{ and } $\theta(1)=\alpha \}$. Then $\text{H}_n(A)$ is a normal subgroup of $\GL_{n}(A)$.

	\item Let $I\subset A$ be an ideal and let ``bar" denote going modulo $I$. We define $\Um_{n}(A,I):=\{v\in \Um_{n}(A):\ol v=e_1\}$ and $\GL_{n}(A,I):=\{\alpha\in \GL_{n}(A):\ol \alpha=\text{Id}\}$. 
	\item Let $P$ be a projective $A$-module. An element $p \in P$ is said to be a unimodular element of $P$ if there exists $\phi \in P^* = \text{Hom}_A(P, A)$ such that $\phi(p) = 1$. The set of all unimodular elements of $P$ is denoted by $\Um(P)$.
	\item  Let $P$ be a projective $A$-module such that $P$ has a unimodular element. We choose $\phi\in P^*$ and $p \in P$ such that $\phi(p)=0$. We define an endomorphism $\phi_p$ as the composite $\phi_p:P\to A\to P$, where $A\to P$ is the map sending $1\to p.$ Then by a \textit{transvection} we mean an automorphism of $P$, of the form $1+\phi_p$, where either $\phi\in \Um(P^*)$ or $p\in \Um(P)$. By $\E(P)$  we denote the subgroup of $\Aut(P)$ generated by all transvections.
	
\end{enumerate}

\ed

We begin by considering the following observation for a graded domain. This simple proposition plays a crucial role in the article, and therefore, we provide the proof.

\bp\label{dogr}
Let $R=\bigoplus_{i\ge 0} R_i$ be a graded domain of dimension $d$. Let $S\subset R_0$ be a multiplicative set such that $S\cap\mm\not=\emptyset$, for any maximal ideal $\mm\in \Spec(R_0)$. Then the graded domain $S^{-1}R$ does not have a graded maximal ideal $S^{-1}\MM$ such that $\MM$ is a maximal ideal in $R$. As a consequence, we get $\dim(S^{-1}R)<d$. 
\ep

\proof Suppose that $\dim(R_0)=n$. We give the proof by induction on $n$. First, we note that if $n=0$, then $R_0$ is a field. In this case, the statement is vacuously true.

Now, we consider the case where $n\geq 1$. If there does not exist such an $S$, then the statement is again vacuously true. Therefore, we assume that such an $S$ exists. Contrarily, we assume the existence of a graded maximal ideal, denoted as $S^{-1}\MM$, in $S^{-1}R$, where $\MM\in \Spec(R)$ is a maximal ideal. There are two possibilities: either $\MM$ is a graded maximal ideal or $\MM$ is a maximal ideal but not a graded ideal. If $\MM$ is a graded maximal ideal, then it can be expressed as $\mm\oplus R_+$, where $\mm$ is a maximal ideal in $R_0$. Since $S\cap \mm\neq\emptyset$, it implies that $S\cap \MM\neq\emptyset$. However, this leads to a contradiction.

Now, we assume that $\MM$ is not a graded ideal. Since $S^{-1}\MM$ is a graded maximal ideal, it can be expressed as $\mm'\oplus S^{-1}R_+$, where $\mm'$ is a maximal ideal in $S^{-1}R_0$. In particular, as $R$ is Noetherian, there exists $s\in S$ such that $sR_+\subset \MM$. Because of $S\cap\MM=\emptyset$, it follows that $R_+\subset \MM$.

We claim that $\mm_0:=\MM\cap R_0\neq{0}$ is a non-zero prime ideal in $R_0$. We observe that, it is enough to show that $\mm_0\not=0$. Contrary, let us assume that $\mm_0=0$. Consider an element $f\in \MM$. We can write $f=f_0+f_1$, where $f_0\in R_0$ and $f_1\in R_+$. Since $f_1\in R_+\subset \MM$, we have $f_0\in \mm_0={0}$. This implies $\MM=R_+$. As a result, we get $R_0\cong R/R_+\cong R/\MM\cong k$, where $k$ is a field. However, this is not possible as $\dim(R_0)=n\geq 1$. Hence, we establish that $\mm_0\neq{0}$.


Let "bar" denote going modulo $\mm_0$. Note that if $\mm_0 \cap S\neq\emptyset$, then we have $S\cap \MM\neq\emptyset$, which contradicts our assumption on the existence of such a maximal ideal. Hence, without loss of generality, we may assume that $\mm_0 \cap S=\emptyset$. Since $\mm_0\not=0$, we have $\dim(\overline{R_0})\leq n-1$, and $\overline{R}$ is a graded domain with dimension $\leq d-1$.

Let $\eta\in \Spec(\overline{R_0})$ be a maximal ideal in $\overline{R_0}$. Since  $\eta+\mm_0$ is also a maximal ideal in $R_0$, according to our induction hypothesis, we have $S\cap \langle \eta+\mm_0\rangle \neq\emptyset$. Implying that $\overline{S}\cap \eta \neq\emptyset$. Furthermore, we observe that $\overline{\MM}$ is a maximal ideal in $\overline{R}$. Moreover, as $\mm_0\subset R_0$ and $\overline S^{-1}\ol \MM\cong\overline{S^{-1} \MM}$, the ideal $\overline S^{-1}\ol \MM$ is a graded maximal ideal in $\overline S^{-1}{\ol R}$ (recall that $S^{-1}\MM$ is a graded ideal). However, by the induction hypothesis, there does not exist such a maximal ideal in $\overline S^{-1}\ol R$. This completes the induction step.

It remains to show that $\dim(S^{-1}R)<d$. To prove this, we note that for an arbitrary graded ring $B=\bigoplus_{i\geq 0} B_i$, there exists a graded maximal ideal $\MN$ in $B$ such that $\hh(\MN)=\dim(B)$. In $S^{-1}R$, any graded maximal ideal of height $d$ is a localization of a maximal ideal in $R$. However, we have already demonstrated the nonexistence of such a graded maximal ideal of height $d$ in $S^{-1}R$. Therefore, the ring $S^{-1}R$ does not have a graded maximal ideal of height $d$. Consequently, we obtain that $\dim(S^{-1}R)<d$.  \qed

We revisit a well-known homotopy map due to Swan and Weibel. 
\smallskip

\bd 
Let $R=\bigoplus_{i\ge 0} R_i$. We define the Swan-Weibel's homotopy map $\Gamma_{SW}: R \to R[T]$ as follows: for any element $f = a_0 + a_1 + \ldots + a_n \in R$, we define $\Gamma_{SW}(f):=a_0+a_1T+\ldots+a_nT^n \in R[T]$, where $a_i \in R_i$.
\ed

\bl\label{isotopy} Let $R=\bigoplus_{i\ge 0} R_i$. Let $\alpha\in \GL_n(R)$ such that $\ol{\alpha}=\text{Id}$, where ``bar" denotes going modulo the ideal $R_+$. Then there exists an $\theta(T)\in \text{GL}_n(R[T])$ such that $\theta(0)=\text{Id}$ and $\theta(1)=\alpha$. In other words $\alpha\in \text{H}_n(R)$. Moreover, if $e_1\alpha=e_1$, then we may choose such an $\theta(T)$ with the property that $e_1\theta(T)=e_1$. 

\el
\proof Consider the group homomorphism $\widetilde{\Gamma_{SW}}:\GL_n(R)\to \GL_n(R[T])$ induced by $\Gamma_{SW}$ \cite[Definition 2.1]{rbs}. Let us take $\theta(T)=\widetilde{\Gamma_{SW}}(\alpha)\in \GL_{n}(R[T])$. Then it follows that $\theta(0)=\ol{\alpha}=\text{Id}$ and $\theta(1)=\alpha$. Now we assume that $e_1\alpha=e_1$. As $1\in R_0$, we have $\Gamma_{SW}(1)=1$. Hence, we have $e_1\theta(T)=e_1$.\qed


\bl\label{SWE}
Let $R=\bigoplus_{i\ge 0} R_i$. Then the map $\Gamma_{SW}:R\to R[T]$ will induce a group homomorphism  $\widetilde{\Gamma_{SW}}:\E_n(R)\to \E_n(R[T])$.
\el
\proof First, we observe that $\Gamma_{SW}$ will induce a group homomorphism $\widetilde{\Gamma_{SW}}:\E_n(R)\to \GL_{n}(R[T])$ for details we refer to \cite[Definition 2.1]{rbs}. Therefore, it is enough to show that $\widetilde{\Gamma_{SW}}(\E_n(R))\subset \E_n(R[T])$. Let $E_{ij}(f)\in \E_n(R)$ be an elementary matrix whose only non-zero non-diagonal entry is $f$ at the position $(i,j)$, where $i\not=j$. We write $f(T)=\Gamma_{SW}(f)$. Then we note that $\widetilde{\Gamma_{SW}}(E_{ij}(f))=E_{ij}(f(T))\in \E_n(R[T])$. Moreover, since $\widetilde{\Gamma_{SW}}$ is a group homomorphism and any element of $\E_{n}(R)$ can be written as a finite product of elements of the form $E_{ij}(f)$, it follows that $\widetilde{\Gamma_{SW}}(\E_n(R))\subset \E_n(R[T])$. \qed

The next lemma is known as one of the variants of Quillen-Suslin's Local-Global Principle and must be well-known. However, we could not find any suitable reference for the exact version required in this article. The closest reference we have found is \cite[Theorem 3.8]{rbs}. Therefore, we provide the proof, which is straightforward using the homotopy map $\Gamma_{SW}$ and Suslin's Local-Global Principle \cite[Lemma 3.5]{AASSC}.

\bl\label{LG}
Let $R=\bigoplus_{i\ge 0} R_i$ and let ``bar" denote going modulo the ideal $R_+$. Let $s,t\in R_0$ be two co-maximal elements and let $\alpha\in \GL_n(R)$ such that (i) $\ol{\alpha}=\text{Id}$, (ii) $\alpha_s\in \E_{n}(R_s)$ and (iii) $\alpha_t\in \E_{n}(R_t)$, where $n\ge 3$. Then $\alpha\in \E_{n}(R)$.
\el
\proof We take $\theta(T)=\widetilde{\Gamma_{SW}}(\alpha)\in \GL_{n}(R[T])$, where $\widetilde{\Gamma_{SW}}:\GL_n(R)\to \GL_n(R[T])$ is induced by $\Gamma_{SW}$. Then we note that $\theta(0)=\text{Id}$. Moreover, it follows from Lemma \ref{SWE} that $(\theta(T))_s\in \E_{n}(R_s[T])$ and $(\theta(T))_t\in \E_{n}(R_t[T])$. Applying \cite[Lemma 3.5]{AASSC} we obtain that $\theta(T)\in \E_{n}(R[T])$. Therefore, we get $\alpha=\theta(1)\in \E_{n}(R)$.\qed

The following lemma is due to Plumstead, which is an adaptation of \cite[Example 4]{P}, tailored to our requirements. Here, we point out that the following version has an additional conclusion compared to the version given in \cite[Lemma 2.4]{sb1}, and this conclusion is crucially used in Lemma \ref{ml}. However, the same proof works here as well. Hence, we omit the proof to avoid repeating the same argument.

\bl\label{plgd}
Let $A$ be a ring of dimension $d$, and let $s$ be a non-zero divisor in $A$ such that $\dim(A_s)\le d-1$. Then there exists a generalized dimension function $\delta:\Spec(A)\to \mathbb{N}$ such that $\delta(\p)\le d-1$ for all $\p\in \Spec(A)$. Furthermore, we can choose $\delta$ such that $\delta(\p)=\dim(A/\p)$ for all  $\p\owns s$.

\el

The next theorem is derived from a pivotal result due to Eisenbud and Evans \cite{EE}. This has been used extensively throughout the article. This version is recollected from \cite[Eisenbud-Evans Theorem]{P}.
\bt\label{eept}
Let $A$ be a ring, and let $\CP\subset\Spec(A)$ be a subset. Consider a generalized dimension function $\delta:\CP \to \mathbb{N}$. Let $M$ be an $A$-module satisfying $\mu_{\p}(M)\ge 1+\delta(\p)$ for all $\p\in \CP$, where $\mu_{\p}(M)$ is the minimal number of generators of $M_{\p}$. For a basic element $(r,m)\in A\oplus M$ on $\CP$, there exists an element $m'\in M$ such that $m+rm'$ is also a basic element on $\CP$.
\et

As a consequence of Theorem \ref{eept}, we have the following result, whose proof can be found in \cite[Corollary 2.13]{BRS}.

\bc\label{ee}

Let $A$ be a ring and $P$ be a projective $A$-module of
rank $n$. Let $(\alpha,a)\in P^*\oplus A$. Then there exists an element $\beta\in P^*$ such 
that $\hh(I_a)\ge n$, where $I:= (\alpha+a\beta)(P)$. 
In particular, if the ideal $\langle \alpha(P),a\rangle $ has height $\ge n$, then $\hh(I)\ge n$. Further, if $\langle \alpha(P),a \rangle$ is an ideal of height $\ge n$ and $I$ is a proper
ideal of $A$, then $\hh(I)=n$.
\ec 

We conclude this section with Quillen’s famous splitting lemma. The proof is essentially contained in \cite[Theorem 1, paragraph 2]{Q} (see also \cite[Lemma 2.9]{MPHIL}).

\bl\label{QSL}
Let $A$ be a ring, and let $s, t \in A$ satisfying $\langle s \rangle + \langle t \rangle = A$. Let $\alpha \in \text{H}_{n}(A_{st})$. Then there exist $\psi_1\in \GL_{n}(A_t)$ and $\psi_2\in \GL_{n}(A_s)$ such that $\alpha=(\psi_1)_s(\psi_2)_t$.
\el

\section{Unimodular rows}\label{ur}
This section is devoted to establishing that any unimodular row of length $d+1$ over a graded domain of dimension $d\ge 1$ can be completed to the first row of an invertible matrix, which is homotopic to the identity matrix. We begin with an easy consequence of Lemma \ref{plgd} and Theorem \ref{eept}.
\bl\label{dsr}
Let $A$ be a ring of dimension $d\ge 1$. Assume that, there exists a non-zero divisor $s\in A$ such that $\dim(A_s)<\dim(A)$. Then $\text{sr}(A)\le d$.
\el
\proof Let $v=(v_1,\ldots,v_{d+1})\in \Um_{d+1}(A)$. Applying Lemma \ref{plgd} we get a generalized dimension function $\delta:\Spec(A)\to \mathbb{N}$ such that $\delta(\p)\le d-1$ for all $\p\in \Spec(A)$. We note that $v$ is a basic element of the free module $A^{d+1}$. Then applying Theorem \ref{eept} (taking $M=A^d$) we obtain a basic element $w=(v_1+\lambda_1v_{d+1},\ldots,v_d+\lambda_dv_{d+1})$ of $A^d$, for some $\lambda_i\in A$. Now since $A^d$ is a free (in particular, a projective) module, every basic element is a unimodular row. This concludes the proof. \qed

The next proposition is similar to the well-known Quillen's splitting lemma \cite[Theorem 1]{Q}. Here we reproduce it in our setup with an added conclusion, which is crucial for this article.

\bp\label{splitting}
Let $R=\bigoplus_{i\ge 0} R_i$ and $s,t\in R_0$ such that $\langle s\rangle +\langle t \rangle =R_0$. Let $\eta\in \GL_n(R_{st},(R_+)_{st})$ such that $e_1\eta=e_1$. Then there exist $\eta_1\in \GL_n(R_s,(R_+)_{s})$ and $\eta_2\in \GL_n(R_t,(R_+)_{t})$ such that 
\begin{enumerate}
	\item  $\eta=(\eta_1)_t(\eta_2)_s$,
	\item $e_1\eta_i=e_1$, for $i=1,2$.
\end{enumerate}
\ep
\proof Let ``bar" denote going modulo the ideal $R_+$. We define $\chi(X):= \widetilde{\Gamma_{SW}}(\eta)$. Using Lemma \ref{isotopy} we obtain the following. $$\chi(X)\in  \GL_{n}(R_{st}[X]) \text{ such that } \chi(0)=\ol\eta= \text{Id} \text{ and }e_1\chi(X)=e_1$$
We claim that there exist $N\in \mathbb N$ and $\lambda\in R_0$  such that if we take $g = \lambda s^N$, then the following holds.
$$\chi(X)\chi(gX)^{-1}\in \GL_{n}(R_s[X])\text{ and }\chi(gX)\in \GL_{n}(R_t[X])$$First we prove our claim. To prove this we follow the argument given in \cite[Lemma 2.9]{MPHIL}. Since $\chi(0)=\text{Id}$, by \cite[Lemma 2.8]{MPHIL} there exists large enough $N_1\in \mathbb N$ such that for all $i\ge N_1$ and for all $\lambda \in R_0$, we have $\chi(\lambda s^iX)\in \GL_{n}(R_t[X])$.

Let us consider two variables $T$ and $Y$. We define $\delta(T,X,Y):=\chi((T+Y)X)\chi(TX)^{-1}$. Then $\delta(T,X,Y)\in \GL_{n}(R_{st}[T,X,Y])$ such that $\delta(T,X,0)=\delta(T,0,Y)=\text{Id}$. Therefore, again applying \cite[Lemma 2.8]{MPHIL} we can find $N_2\in \mathbb N$ such that for all $j\ge N_2$ and for all $\mu\in R_0$, we have ${\delta}(T,X, t^j\mu Y)\in \GL_{n}(R_s[T,X,Y])$.

Let us choose $N=\max\{N_1,N_2\}$. As $\langle s \rangle +\langle t\rangle =R_0$, there exist $\lambda,\mu \in R_0$ such that $\lambda s^N+\mu t^N=1$. Let us take $g=\lambda s^N$. Now we write $\chi(X)=\chi(X)\chi(gX)^{-1}\chi(gX)$. By our choice of $N$, we have $\chi(gX)=\chi(\lambda s^NX)\in \GL_{n}(R_t[X])$. Now one may observe the following. $$\chi(X)\chi(gX)^{-1}=\chi((g+\mu t^N)X)\chi(gX)^{-1}=\delta(g, X, \mu t^N)\in \GL_n(R_s[X])$$This proves our claim.

Since $e_1\chi(X)=e_1$, we further obtain that $e_1\chi(gX)=e_1$ and $e_1\chi(X)\chi(gX)^{-1}=e_1$. Let us define $\eta_1:=\chi(1)\chi(g)^{-1}$ and $\eta_2:=\chi(g)$. We observe the matrix $\chi(X)$ has the property that $\ol{\chi(a)}=\text{Id}$ for any $a\in R_0$. Since $g\in R_0$, this further imply that $\ol \eta_i=\text{Id}$, for $i=1,2$. Therefore, we get the following.
\begin{enumerate}[\quad \quad (1)]
	\item  $\eta=(\eta_1)_t(\eta_2)_s$;
	\item   $\eta_1\in \GL_{n}(R_s,(R_+)_{s})$;
	\item $\eta_2\in \GL_{n}(R_t,(R_+)_{t})$;
	\item  $e_1\eta_i=e_1$, for $i=1,2$.
\end{enumerate}This concludes the proof.\qed 

\smallskip

The next lemma concerns the patching of two invertible matrices in a graded ring.

\bl\label{patching}
Let $R=\bigoplus_{i\ge 0} R_i$ and $s,t\in R_0$ such that $\langle s\rangle +\langle t \rangle =R_0$. Let $v\in \Um_n(R,R_+)$. Assume that, there exist $\alpha_1\in \GL_{n}(R_s,(R_+)_s)$ and $\alpha_2\in \GL_n(R_t,(R_+)_t)$ such that $v\alpha_i=e_1$, for $i=1,2$. Then there exists an $\alpha\in \GL_n(R,R_+)$ such that $v\alpha=e_1$.
\el

\proof Let ``bar'' denote going modulo $R_+$. Let us define $\eta:=(\alpha_1)_t^{-1}(\alpha_2)_s\in \GL_{n}(R_{st})$. Then we note that $\ol\eta=\text{Id}$ and $e_1\eta=e_1$. Applying Proposition \ref{splitting} there exist $\eta_1\in \GL_n(R_s,(R_+)_{s})$ and $\eta_2\in\GL_n(R_t,(R_+)_{t})$ such that 
\begin{enumerate}[\quad \quad (a)]
	\item  $\eta=(\eta_1)_t(\eta_2)_s$,
	\item $e_1\eta_i=e_1$, for $i=1,2$.
\end{enumerate}
We now define $\sigma_1:=\alpha_1\eta_1\in \GL_{n}(R_s,(R_+)_s)$ and $\sigma_2:=\alpha_2\eta_2^{-1}\in \GL_{n}(R_t,(R_+)_t)$. Here we notice that $v\sigma_i=e_1$ $(i=1,2)$. Because of $\eta=(\alpha_1)_t^{-1}(\alpha_2)_s =(\eta_1)_t(\eta_2)_s$, we have $(\sigma_1)_t=(\sigma_2)_s$. Therefore, by \cite[Proposition 2.2, page no 211]{Lam} there exists a unique $\alpha\in \GL_{n}(R,R_+)$ such that $\alpha_s=\sigma_1$ and $\alpha_t=\sigma_2$. Furthermore, the matrix $\alpha$ takes $v$ to $e_1$ as it is true locally. \qed

\smallskip

\notation Let $A$ be a ring. 
\begin{enumerate}[(i)]
	\item Let $G\subset \GL_{n}(A)$ be a subgroup. For any $u,v\in \Um_n(A)$, we define $ u \sim_{G} v$ if there exists an $\epsilon\in G$  such that $u\epsilon=v$. We denote the set $\{v\in \Um_n(A): v \sim_{G} e_1\}$ by the notation $e_1G$.
	
	\item The Jacobson radical of $A$ is denoted by $\text{Jac}(A)$.
\end{enumerate}

\smallskip

\rmk\label{rpvul} One may wonder whether it is possible to improve Quillen's splitting lemma (for an arbitrary ring) in such a way that both the splitting matrices fix the canonical vector $e_1$. Unfortunately, achieving such an improvement, as claimed in \cite[Lemma 3.7]{RIP}, is not feasible. To illustrate this, here we argue as follows: consider a ring $A$ of dimension $d\ge 2$. We show that such an improvement of Quillen's splitting lemma will ultimately lead to the conclusion that $\Um_{d+1}(A)=e_1\text{SL}_{d+1}(A)$. However, this is not true as this discrepancy is illustrated by the well-known example of the projective module corresponding to the tangent bundle of an even-dimensional real sphere. To establish the mentioned implication we choose a $v\in \Um_{d+1}(A)$. Then one can always find a non-zero divisor $s\in R$ such that $v\alpha_1= e_1$, for some $\alpha_1\in {\E_{d+1}(A_s)}$. As $s\in \text{Jac}(A_{1+\langle s \rangle })$ is a non-zero divisor, it is not difficult to establish that $v\alpha_2= e_1$, for some $\alpha_2\in {\E_{d+1}(A_t)}$ and $t\in 1+\langle s \rangle$. Now, if the elementary matrix $\eta=(\alpha_1)_t^{-1}(\alpha_2)_s$ splits in such a way that each of its splitting matrices fixes $e_1$, then applying the arguments given in Lemma \ref{patching} one can obtain an $\alpha\in \SL_{d+1}(A)$ such that $v\alpha=e_1$. This, in particular shows that $\Um_{d+1}(A)=e_1\text{SL}_{d+1}(A)$.

\rmk One may observe that in \cite{RIP}, to prove one of their main results, Theorem 3.8, Lemma 3.7 (of the same article) plays a very crucial role. Moreover, in the same article, Theorem 3.8 plays a significant role in establishing results in Section 4. However, a completely more general result compared to \cite[Theorem 3.8]{RIP} has been proved independently in \cite[Theorem 4.5]{sb1}.

\bt\label{umch}
Let $R=\bigoplus_{i\ge 0} R_i$ be a graded domain of dimension $d\ge 2$ such that $\dim(R_0)\ge 1$. Then for any $v\in \Um_{d+1}(R,R_+)$ there exists an $\alpha\in \GL_{d+1}(R,R_+)$ such that $v\alpha=e_1$. As a consequence
$$ \Um_{d+1}(R)=e_1\text{H}_{d+1}(R).$$
\et
\proof Let $v\in \Um_{d+1}(R)$ and let ``bar" denote going modulo $R_{+}$. As $R$ has a non-trivial grading the ideal $R_+\not=0$. In particular, we get $\hh(R_+)\ge 1$. As $\dim(\ol R)<d$, one can use Prime avoidance lemma to prove that $\text{sr}(\ol R)\le d$ (cf. \cite[Theorem 3.5, $\S$ 3, page no 239]{Bass68}). Hence there exists $\ol\kappa\in \E_{d+1}(\ol R)$ such that $\ol v\ol \kappa=\ol e_1$. Since the canonical map $\E_{d+1}(R)\twoheadrightarrow \E_{d+1}(\ol R)$ is surjective, there exists a lift $\kappa\in \E_{d+1}(R)$ of $\ol{\kappa}$. Altering $v$ by $v\kappa$ one may further assume that $\ol v=\ol e_1$. Now, if there exists an $\alpha\in \GL_{d+1}(R,R_+)$ such that $v\alpha=e_1$, then it follows from Lemma \ref{isotopy} that $ \Um_{d+1}(R)=e_1\text{H}_{d+1}(R).$ Hence, to prove the theorem it is enough to find such an $\alpha$. In the remaining part of the proof we find such an $\alpha$.

Let $S=R_0\setminus\{0\}$. Applying Proposition \ref{dogr} we get $\dim(S^{-1}R)\le d-1$. Again using \cite[Theorem 3.5, $\S$ 3, page no 239]{Bass68} we obtain that $\text{sr}(S^{-1}R)\le d$. Therefore, we can find an $s\in S$ such that $v\sim_{\E_{d+1}(R_s)} e_1$. Let $\alpha_1\in \E_{d+1}(R_s)$ be such that $v\alpha_1=e_1$. Furthermore, we may replace $\alpha_1$ by $\alpha_1 \ol\alpha_1^{-1}$ and assume that $\ol\alpha_1=\text{Id}.$

Let $\MT=\{1+sr:r\in R_0\}$ and $B=\MT^{-1}R$. Since $\MT\subset R_0$, the ring $B$ is also a graded ring. Moreover, we note that $s\in \text{Jac}(\MT^{-1}R_0)$. Hence, applying Proposition \ref{dogr} it follows that $\dim(B_s)\le d-1.$ Therefore, by Lemma \ref{dsr} we get $\text{sr}(B)\le d$. Thus, there exists an $\alpha_2\in \E_{d+1}(B)$ such that $v\alpha_2=e_1$. As again we may replace $\alpha_2$ by $\alpha_2 \ol\alpha_2^{-1}$ and further assume that $\ol\alpha_2=\text{Id}.$ We can find $t\in \MT$ such that $\alpha_2\in \E_{d+1}(R_t)$.

Now applying Lemma \ref{patching} we can find an $\alpha\in \GL_{d+1}(R,R_+)$ such that $v\alpha=e_1$. This completes the proof.\qed

\smallskip

\rmk\label{umchrmk} Let $A$ be a regular ring of essentially finite type over a field. Then using \cite[Theorem 3.3]{TV} it follows that $\text{H}_{n+1}(A)=\E_{n+1}(A)$, for all $n\ge 2$. Hence, in Theorem \ref{umch}, additionally if we assume that $R$ is a regular ring of essentially finite type over a field, then we get $\Um_{d+1}(R)=e_1\E_{d+1}(R)$. However, we do not know whether the regularity of $R$ is actually necessary.
\smallskip

\rmk \label{Rmk1}One can remove the hypothesis that $\dim(R_0)\ge 1$ in Theorem \ref{umch} at the expense of the hypothesis that $\frac{1}{d!}\in R$ by utilizing the Swan-Weibel's homotopy map and applying \cite[Corollary 2.5]{Rao}. In fact, the same yields the following: let $R=\bigoplus_{i\ge 0} R_i$ be a graded ring of dimension $d$ such that $\frac{1}{d!}\in R$. Then $\Um_{d+1}(R)=e_1\text{H}_{d+1}(R)$. It is worth noting that the removal of the hypothesis ``$\frac{1}{d!}\in R$" from a cancellation problem is highly non-trivial (cf. \cite{AAS} and \cite{FRS}).


In the remaining part of the section, we extend Theorem \ref{umch} over various over-rings of the graded rings considered in Theorem \ref{umch}.

\bc\label{ec1}
Let $R$ and $d$ be as in Theorem \ref{umch} and let $A=\MS^{-1}R$, where $\MS\subset R$ is a multiplicative set. Then $$\Um_{d+1}(A)=e_1\text{H}_{d+1}(A).$$
\ec
\proof First we comment that, since $\E_{d+1}(A)\subset \text{H}_{d+1}(A)$, the only non-trivial case is when $\dim(A)=d$. Hence, without loss of generality, we assume that $\dim(A)=d$. Let us choose $v\in \Um_{d+1}(A)$. Then there exists a non-zero $x\in R$ such that (i) $v\in \Um_{d+1}(R_x)$ and (ii) $\dim(R_x)=d$. Now it follows from \cite[Lemma 4.4]{sb1} that there exists $u\in \Um_{d+1}(R)$ such that $v\sim_{\E_{d+1}(R_x)}u$. Applying Theorem \ref{umch} we obtain that $u\sim_{\text{H}_{d+1}(R)} e_1$. Since $\sim_{\text{H}_{d+1}(R_x)}$ is transitive, the proof concludes. \qed

\bc\label{ec2}

Let $R$ and $d$ be as in Theorem \ref{umch}. Additionally, we assume that $R$ is an affine domain over a field. Then 
$$\Um_{d+1}(R[X_1,\ldots,X_n])=e_1\SL_{d+1}(R[X_1,\ldots,X_n]).$$
\ec
\proof We use Quillen Induction on $n$ to prove the theorem. For $n=0$ this follows from Theorem \ref{umch}. Now let us assume that $n>0$. Let $v\in \Um_{d+1}(R[X_1,\ldots,X_n])$. We note that $R[X_1]=\bigoplus_{i\ge 0} R_i[X_1]$, where $R_i[X_1]=\{\sum_{j=1}^n a_jX_1^j: a_j\in R_i \text{ and } n\in \mathbb N \}$ and the sum is defined in the obvious way. As there exists a canonical surjection $R\twoheadrightarrow R_0$, the ring $R_0$ is also an affine domain over the same field, say $k$. Let us take $S=k[X_1]\setminus \{0\}\subset R_0[X_1]$. Therefore, we have $\dim(S^{-1}R_0[X_1])=\dim(R_0)$. As any maximal ideal $\mm$ of $R_0[X_1]$ is of height equal to $\dim(R_0)+1$, we have $S\cap \mm\not=\emptyset$. We take $B=S^{-1}R[X_1]$. Then it follows from Proposition \ref{dogr} that $\dim(B)\le d$. In particular, since $\dim(S^{-1}R_0[X_1])=\dim(R_0)$ we have $\dim(B)=d$. Then $B=\bigoplus_{i\ge 0} S^{-1}R_i[X_1](=\bigoplus_{i\ge 0} B_i$ say$)$ is also a graded affine domain over the field $k(X_1)$ of dimension $d$ such that $\dim(B_0)=\dim(R_0)\ge 1$. Applying induction hypothesis on $B[X_2,\ldots,X_n]$ we can find a monic polynomial $f\in S$ such that $$v\sim_{\SL_{d+1}(D[X_1]_f)} e_1,$$ where $D=R[X_2,\ldots,X_n]$. Then by Affine Horrocks Theorem \cite[Theorem 3]{Q} the result follows.\qed

We end this section with a theorem on the injective stability of $\text{K}_1(R)$, where $R$ is a graded domain. The proof is a straightforward consequence of Lemmas \ref{LG} and \ref{dsr}. Before that, we restate a stability theorem due to Vaser{\v{s}}te{\u{\i}}n to suit our needs. One can find the proof in \cite[Theorem 3.2]{Va1}.

\bt[Vaser{\v{s}}te{\u{\i}}n]\label{VST}
Let $A$ be a ring with $\text{sr}(A) \leq m$. Then, for any $n \geq m + 1$, we have $\SL_n(A) \cap \E(A) = \E_n(A)$.
\et

	\bt\label{ist}
Let $R$ and $d$ be as in Theorem \ref{umch}. Additionally, for $n\ge 1$ we further assume that $R$ is an affine domain over a field. Then $$\SL_{d+1}(R[X_1,\ldots,X_n])\cap \E(R[X_1,\ldots,X_n])=\E_{d+1}(R[X_1,\ldots,X_n]).$$
\et

\proof We again apply Quillen Induction on $n$ to prove the theorem. We give the proof in cases.\\
\textbf{Case - 1.} Let us assume that $n=0$. Let $\alpha\in \SL_{d+1}(R)\cap \E(R)$ and let ``bar'' denote going modulo the ideal $R_+$. Since $\ol\alpha\in \SL_{d+1}(R_0)$ and $R_0\subset R$, we may treat $\ol\alpha$ as an element of $\SL_{d+1}(R)$. Moreover, we observe that as $\hh(R_+)\ge 1$ we have $\dim(R_0)<\dim(R)$. In other words, we get $\text{sr}(R_0)\le d$. Hence applying Theorem \ref{VST} we obtain that (1) $\alpha\in \SL_{d+1}(R) \cap \E_{d+2}(R)$ and (2) $\ol \alpha\in \E_{d+1}(R_0)\subset \E_{d+1}(R)$. We take $\beta=\alpha\ol\alpha^{-1}$. Then we note that $\beta \in \SL_{d+1}(R) \cap \E_{d+2}(R).$ Let us consider $\MT=R_0\setminus \{0\}$. Then applying Proposition \ref{dogr} we get $\dim(\MT^{-1}R)\le d-1$. Hence, again using Theorem \ref{VST} on $\MT^{-1}R$ we obtain that 
$(\beta)_{\MT}\in \E_{d+1}(\MT^{-1}R)$. There exists an $s\in \MT$ such that $\beta_s\in \E_{d+1}(R_s)$.

Let $S=\{1+sr:r\in R_0\}$ and let $B=S^{-1}R$. Then again by Proposition \ref{dogr} we obtain that $\dim(B_s)\le d-1$. Hence, using Lemma \ref{dsr} it follows that $\text{sr}(B)\le d$. We again apply Theorem \ref{VST} to obtain that $(\beta)_S\in \E_{d+1}(B)$. We choose $t\in S$ such that $\beta_t\in \E_{d+1}(R_t)$. Now it follows from Lemma \ref{LG} that $\beta\in \E_{d+1}(R)$. Because of $\ol\alpha\in \E_{d+1}(R)$, we have $\alpha\in \E_{d+1}(R)$. This concludes the proof for $n=0$.\\
\textbf{Case - 2.} Now let us assume that $n>0$. Then applying Quillen Induction on $n$ as described in Corollary \ref{ec2} and using \cite[Corollary 5.7]{AASSC} suitably one may conclude the proof.\qed


\section{Main theorems}\label{MT}
Let $A$ be a ring and $I\subset A$ be an ideal. We call $I$ is efficiently generated if $\mu(I/I^2)=\mu(I)$. This section is devoted to studying the efficient generation problem for top height ideals in a ring. Before presenting the main theorems, we need some preparation. We begin this section with the following lemma, which is a consequence of Lemma \ref{plgd} and Theorem \ref{eept}.

\bl\label{ml}
Let $A$ be a ring of dimension $d\ge 2$. Assume that, there exists a non-zero divisor $s\in A$ such that $\dim(A_s)<d$. Let $I\subset A$ be an ideal such that $\mu(I/I^2)=\hh(I)=d$. Then any set of generators of $I=\langle f_1,\ldots,f_d\rangle +I^2$ lifts to a set of generators of $I$.
\el

\proof Applying  Lemma \ref{plgd} one may obtain a generalized dimension function $\delta:\Spec(A) \to \mathbb{N}$ such that $\delta(\mq)\le d-1$ for all $ \mq\in \Spec(A)$ and $\delta(\mq)=\dim(A/\mq )$ for all $\mq\owns s$. Let $\p\in \Spec(A)$. Suppose that, we have $I\subset \p$. As $\hh(I)=d$ we must have $\hh(\p)=d$. Since $\dim(A_s)<d$, the element $s$ is in $\p$. This implies $\delta(\p)=\dim(A/\p)=0$. As $A_{\p}$ is a local ring we have $\mu(IA_{\p}/I^2A_{\p})=\mu(IA_{\p})=d$. Therefore, we obtain that $\mu(IA_{\p})+\delta(\p)\le d$. Now if $I\not \subset \p$, then $\mu(IA_{\p})=1$. Thus also in this case we have $\mu(IA_{\p})+\delta(\p)\le d$. In particular, we get $\sup\{\mu(IA_{\p})+\delta(\p):\p\in \Spec(A)\}\le d$. Hence, one may apply \cite[Theorem 0]{P} to find $e_i\in I^2$ $(i=1,\ldots,d)$ such that $I=\langle l_1,\ldots,l_d\rangle $, where $l_i=f_i+e_i$ . This completes the proof.\qed

In the following proposition, we present a general criterion for the efficient generation of a top height ideal in an arbitrary ring. This criterion enables us to identify the essential requirements to apply Mohan Kumar's fundamental technique presented in \cite{NMK} to solve the efficient generation problem. By doing so, we are able to provide a unified approach in Theorem \ref{eg}.

\bp\label{egc}
Let $A$ be a ring of dimension $d\ge 2$. Let $I\subset A$ be an ideal such that $\mu(I/I^2)=\hh(I)=d$. Suppose that $I=\langle f_1,\ldots,f_d\rangle +I^2$. Moreover, we assume that there exists a non-zero divisor $s\in A$ and a multiplicative set $S\subset \{1+sr:r\in A\}$ such that the following hold.
\begin{enumerate}
	\item   $IA_s=\langle f_1,\ldots,f_d \rangle A_s+I^2A_s$, has a lift to a set of generators of $IA_s$,
	\item  $\dim(S^{-1}A_s)<d$, and
	\item $\Um_d(S^{-1}A_s)=e_1\text{H}_d(S^{-1}A_s)$.
\end{enumerate}
Then there exist $F_i\in A$ such that $I=\langle F_1,\ldots,F_d\rangle $, with $f_i-F_i\in I^2$. 
\ep
\proof Let $g_i\in IA_s$ be such that $IA_s=\langle g_1,\ldots,g_d\rangle A_s$ with $f_i-g_i\in I^2A_s$. We observe that, if $\dim(A_s)<d$, then applying Lemma \ref{ml} the proof follows. Hence, we assume that $\dim(A_s)=d$. Suppose that $s\not\in \sqrt{I}$, then one may note that $d=\dim(A_s)\ge \hh(IA_s)\ge \hh(I)=d$. Now we wish to apply \cite[Lemma 3.1]{BSB} which is a modification \cite[Lemma 5.6]{SMBB3}. First, we briefly discuss the conclusion of \cite[Lemma 3.1]{BSB} without using the language of the Euler class group. It states that if any ideal $J \subset A$ of height $d$, with a set of $d$-generators of $J/J^2$, say $\omega_J$, and $t \in A$, satisfies all the hypotheses given there, then one can find another ideal $K \subset A$ of height $d$ and a set of $d$-generators of $K/K^2$, say $\omega_K$, such that (a) $t \in \sqrt{K}$, and (b) $\omega_J$ lifts to a set of generators of $J$ if and only if $\omega_K$ lifts to a set of generators of $K$ (the actual statement is more general there). Note that with (1) all the hypotheses of \cite[Lemma 3.1]{BSB} are satisfied. Hence applying the same without loss of generality we may assume that $s\in \sqrt{I}$.

 Let us take $B=S^{-1}A$. Since $\dim(B_s)< d$, using Lemma \ref{ml} we can lift $f_i$'s to a set of generators of $IB$. In particular, we get $l_i\in IB$ such that $IB=\langle l_1,\ldots,l_d\rangle $ and $l_i-f_i\in I^2$, for $i=1,\ldots,d$.  

Since $s\in \sqrt{I}$, the row vectors $(g_1,\ldots,g_d) \text{ and } (l_1,\ldots,l_d)$ are in $ \Um_{d}(B_s)$. Hence, by hypothesis (3) there exists an $\epsilon \in \text{H}_{d}(B_s)$ such that $(g_1,\ldots,g_d)\epsilon= (l_1,\ldots,l_d)$. As $\epsilon \in \text{H}_{d}(B_s)$ there exists a $\theta(T)\in \GL_{n}(B_s[T])$ such that $\theta(0)=\text{Id}$  and $\theta(1)=\epsilon.$ Since $A$ is a Noetherian ring and there are only finitely many $g_i$ and $l_i$, we can find $t\in S$ such that
\begin{enumerate}[\quad\quad (1)]
	\item $IA_t=\langle l_1,\ldots,l_d\rangle$ with $f_i-l_i\in I^2A_t$;
	\item $\theta(T)\in \GL_{d}(A_{st}[T])$. 
\end{enumerate}
Applying Quillen's splitting lemma \ref{QSL} we obtain $\epsilon_1\in \GL_{d}(A_s)$ and $\epsilon_2\in \GL_d(A_t)$ such that $\epsilon=(\epsilon_1)_t(\epsilon_2)_s$. Because of $\langle s\rangle +\langle t \rangle =A$, one may apply a standard patching to obtain $F_i\in I$ such that $I=\langle F_1,\ldots,F_d\rangle $ with $F_i-f_i\in I^2$ for $i=1,\ldots,d$.\qed


Now we are ready to prove the main theorem of this section.

\bt\label{eg}
Let $R=\bigoplus_{i\ge 0} R_i$ be a graded domain of dimension $d$ such that $\dim(R_0)\ge 1$. Let $A$ and $n$ be one of the following:
\begin{enumerate}
	\item  $A=R$ and $n=\dim(R)=d\ge 3$.
	\item $A=\MS^{-1}R$, where $\MS\subset R$ is a multiplicative set contained in the set of all non-zero divisors in $R$ such that $\dim(A)=\dim(R)$ and $n=\dim(R)=d\ge 3 $.
	\item $A=B[M]$, where $B$ is a ring of dimension $ \ge 2$ and $M$ is a finitely generated commutative cancellative (not necessarily torsion free) monoid of rank $r\ge 1$. We take $n=\dim(B[M])$.
\end{enumerate}
Let $I\subset A$ be an ideal such that $\mu(I/I^2)=\hh(I)=n$. Then any set of generators of $I=\langle f_1,\ldots,f_n\rangle +I^2$ lifts to a set of generators of $I$.
\et

\proof We will show that for each of the above rings all the hypotheses of Proposition \ref{egc} are satisfied. We handle these three rings separately in the following cases.\\
\textbf{Case - 1.} In this case, we assume that $A=R$. Let us take $\MT=R_0\setminus \{0\}$. Then by Proposition \ref{dogr} the dimension of $\MT^{-1}A$ is strictly smaller than $n$. Hence, applying \cite{NMK} we can lift $f_i$'s to a set of generators of $\MT^{-1}I$. Therefore, there exist $s\in \MT$ and $g_i\in A_s$ with $I_s=\langle g_1,\ldots,g_n \rangle $ such that $f_i-g_i\in I_s^2$ for $i=1,\ldots,n$. Let $S:=\{1+sr:r\in R_0\}\subset \{1+sx:x\in A\}$. Then again applying Proposition \ref{dogr} we have $\dim(S^{-1}A_s)<n$. Moreover, we observe that since $S \subset R_0$ and $s \in R_0$, the ring $S^{-1}A_s$, which is the same as $S^{-1}R_s$, retains the grading induced by $R$. Because of $n=d\ge 3,$ using Theorem \ref{umch} we have $\Um_{n}(S^{-1}R_s)= e_1\text{H}_{n}(S^{-1}R_s)$. Therefore, applying Proposition \ref{egc} we obtain the required lift.\\
\textbf{Case - 2.} In this case we assume that $A=\MS^{-1}R$. Let $\MT$ be as considered in Case - 1. Then as it was shown in the previous case that $\dim(\MT^{-1}R)<n$, which further implies that $\dim(\MT^{-1}A)<n$. Therefore, following the arguments in the previous case, we can find a non-zero divisor $s\in R_0$ and $g_i\in A_s$ such that $I_s=\langle g_1,\ldots,g_n\rangle $, with $f_i-g_i\in I_s^2$. Let us take $S=\{1+sr:r\in R_0\}\subset \{1+sx:x\in A\}$. Then $S^{-1}A_s=\MS^{-1}(S^{-1}R_s)$, where $S^{-1}R_s$ is a positively graded ring of dimension $\le n-1$. Hence, applying Corollary \ref{ec1} we get $\Um_{n}(S^{-1}A_s)=e_1\text{H}_{n}(S^{-1}A_s)$. Now one may apply Proposition \ref{egc} to complete the proof.\\
\textbf{Case - 3.} In this case we take $A=B[M]$. First we note that for a monoid ring $B[M]$ we  have $\dim(B[M])=\dim(B)+\rank(M)$ \cite[Theorem 4.23]{gulbook}. Let $\MT$ be the set of all non-zero divisors in $B$. Then $\dim(\MT^{-1}B[M])=r$. Since $n>\dim(\MT^{-1}B[M])$ by \cite{NMK} we can lift $f_i$'s to a set of generators of $\MT^{-1}I$. Therefore, there exist $s\in \MT$ and $g_i\in B_s[M]$ with $I_s=\langle g_1,\ldots,g_n \rangle $ such that $f_i-g_i\in I_s^2$ for $i=1,\ldots,n$. Let $S=\{1+sr:r\in B\}$ and let $C=S^{-1}B$. Then as $\dim(C_s)<\dim(B)$ we have $\dim(S^{-1}B_s[M])=\dim(C_s[M])<n$. Therefore, using \cite[Theorem 1.1]{gul} we get $\Um_{n}(C_s[M])= e_1\E_{n}(C_s[M])$. Now one may apply Proposition \ref{egc} to complete the proof.\qed

We now provide an example that proves the necessity of the hypothesis $\dim(R_0)>0$ in Theorem \ref{eg} (1). We essentially use the example constructed by Bhatwadekar, Mohan Kumar and Srinivas \cite[Example 6.4]{BR} in which they provided a non-smooth graded domain (with the degree zero subring a field) over which Nori's question has a negative answer. 

\example\label{ex1} Consider the graded domain $B=\frac{\mathbb{C}[X,Y,Z,W]}{\langle X^5+Y^5+Z^5+W^5\rangle }$. By  \cite[Example 6.4]{BR} there exist (1) an ideal $I\subset B[T]$ such that $\mu(I/I^2T)=\hh(I)=3$ and (2) a set of generators $I=\langle f_1,f_2,f_3\rangle +I^2T$, which does not lift to a set of generators of $I$. Let $S=\mathbb{C}[T]\setminus \{0\}$ and $C=S^{-1}B[T]$. Moreover, one may observe that $I$ does not contain a monic polynomial in $T$. As if it did, then by \cite[Theorem 2.1]{SM}, one could lift $I=\langle f_1,f_2,f_3\rangle +I^2T$ to a set of generators for $I$. Then $C$ is a graded domain of dimension $3$ such that the degree zero subring of $C$ is the field $\mathbb{C}(T)$, and $IC$ is an ideal of $C$ of height $3$ such that $IC=\langle f_1,f_2,f_3\rangle C+I^2C$ does not lift to a set of generators of $IC$. As if such a lift exists then by \cite[Theorem 3.10]{MKD1} one can lift $I=\langle f_1,f_2,f_3\rangle +I^2T$ to a set of generators of $I$, which is not true by (2).
\smallskip


The following is an interesting consequence of the previous theorem. For monoid rings, this is an improvement of \cite[Theorem 3.4]{mkmm}.

\bc\label{eue}
Let $A$ and $n$ be as in Theorem \ref{eg}. Let $P$ be a projective $A$-module with trivial determinant of rank $n$. Then $P$ has a unimodular element.
\ec

\proof Let us choose $\alpha\in P^*$. If $\alpha(P)=A$, then there is nothing to prove. Hence without loss of generality we may assume that $\alpha(P)\subsetneq A$. We apply Corollary \ref{ee} on the pair $(\alpha,1)$ to obtain an element $\beta\in P^*$ such that $\hh(I)\ge n$, where $I=(\alpha+\beta)(P)$. Now again, if $I=A$, then the result is proved. Therefore, without loss of generality we may assume that $I$ is a proper ideal in $A$. Hence it follows from Corollary \ref{ee} that $\hh(I)=d$. Now the result follows from Theorem \ref{eg} and applying subtraction principle \cite[Corollary 3.5]{SMBB3}.\qed
\smallskip

\example\label{ex2} Here we show that the hypothesis $\dim(R_0)\ge 1$ is also necessary in Corollary \ref{eue}, where $A=\bigoplus_{i\ge 0} R_i$. Let $C,S,B,I$ and $f_i$ be as in Example \ref{ex1}. Recall that the $d$-th Euler class group $E^d(D[T])$ and the weak Euler class group $E_0^d(D[T])$ as defined in \cite{MKD1}, where $D$ is a ring of dimension $d\ge 3$ such that $\mathbb Q\subset D$. We consider $(I,\omega_I)\in E^d(B[T])$, where $\omega_I$ is the local orientation induced by $I=\langle f_1, f_2, f_3\rangle +I^2$. Applying \cite[Theorem 2.7]{SMBRS2} we can find a projective $B[T]$-module $P$ (with trivial determinant) of rank $3$ and a surjection $\theta:P\twoheadrightarrow I$. We claim that $S^{-1}P$ does not have a unimodular element. First, we note that if $S^{-1}P$ has a unimodular element, then there exists $f\in \mathbb C[T]\setminus \{0\}$ such that $P_f$ has a unimodular element. But then it follows from \cite[Theorem 3.4]{BRS} that $P$ has a unimodular element. Hence, to prove our claim it is enough to show that $P$ does not have a unimodular element. We fix a trivialization $\chi:\wedge^3 P\iso B[T]$. Then it follows from \cite[Proposition 5.8]{MKD1} and \cite[Theorem 3.4]{SBMKD2} that $E^3(B[T])\cong E^3_0(B[T])$. In particular, this give us $e(P,\chi)=(I,\omega_I)$ in $E^3(B[T])$. Moreover, using \cite[Corollary 4.11]{MKD1} we obtain that $P$ has a unimodular element if and only if $(I,\omega_I)=0$ in $E^3(B[T])$. Now if $(I,\omega_I)=0$ in $E^3(B[T])$, then one may also lift $IC=\langle f_1,f_2,f_3\rangle C+I^2C $ to a set of generators of $IC$. However, as it is shown in Example \ref{ex1} that this is not feasible. Hence, the module $P$ does not have a unimodular element.

In the next theorem we extend Theorem \ref{eg} and Corollary \ref{eue} to the case where the dimension of the graded ring is $2$. 
%

\bt\label{egd2}
Let $R$ be as in Theorem \ref{umch} and $\dim(R)=2$. Let $P$ be a projective $R$-module of rank $2$ with trivial determinant. Suppose $I\subset R$ is an ideal such that $I=\langle f_1, f_2\rangle +I^2$. Then 
\begin{enumerate}[\quad\quad (1)]
	\item $P$ is a free module and
	\item there exist $F_i\in I$ such that $I=\langle F_1,F_2\rangle $, with $F_i-f_i\in I^2$.
\end{enumerate}
\et
\proof We consider $S=R_0\setminus \{0\}$. Then by Lemma \ref{dogr} we get $\dim(S^{-1}R)\le 1$. Since determinant of $P$ is trivial, it follows from \cite{S} that the module $S^{-1}P$ is free. As $P$ is finitely generated module over a Noetherian ring there exists an $s\in S$ such that $P_s$ is a free module. Let us take $T=\{1+sr:r\in R_0\}$ and $B=T^{-1}R$. Then applying Lemma \ref{dogr} we obtain that $\dim(B_s)\le 1$. By Lemma \ref{plgd} and Theorem \ref{eept} we get that, the module $T^{-1}P$ has a unimodular element. Moreover, as determinant of $P$ is trivial, the module $T^{-1}P$ is free. Thus, there exists an element $t\in T$ such that $P_t$ is a free module. Therefore, the Quillen ideal $J(R_0,P)$ of $P$ is $R_0$. In other words, we have $P\cong \frac{P}{PR_+}\otimes R$.  As $\dim(R_0)=1$, again by \cite{S} the $R_0$-module $\frac{P}{PR_+}$ has a unimodular element and hence free (as determinant of $P\cong \frac{P}{PR_+}$ is trivial). Implying that $P$ is free. 

Now we consider $I=\langle f_1,f_2\rangle +I^2$. By \cite{NMKL} there exists $e\in I^2$ such that $I=\langle f_1,f_2,e\rangle $ where $e(1-e)\in \langle f_1,f_2\rangle $. Then $I_e=R_e=\langle 1,0\rangle$ and $I_{1-e}=\langle f_1,f_2\rangle_{1-e}$. Since any unimodular row of length two can be completed to an invertible matrix, using a standard patching argument we obtain a projective $R$-module $Q$ of rank $2$ with trivial determinant and a surjection $\gamma:Q\twoheadrightarrow J$ such that $\gamma$ locally lifts $\{f_1,f_2\}$. Now as $Q$ is free by the previous case, the theorem concludes.\qed


\subsection{Precise obstruction} Let $R = \bigoplus_{i \geq 0} R_i$ be a graded domain of dimension $d \geq 3$, satisfying $\dim(R_0) = 0$ and $\frac{1}{d!} \in R$. In this subsection, we aim to provide a necessary and sufficient condition for the splitting of projective $R$-modules of rank $d$ with a trivial determinant. As the idea used in the proof is well-established due to R. Sridharan \cite{RS}, we will present only a sketch of the proof.

\bt
Let $R=\bigoplus_{i\ge 0} R_i$ be a graded domain of dimension $d\ge 3$ such that $\dim(R_0)= 0$ and $\frac{1}{d!}\in R$. Let $P$ be a projective $R$-modules of rank $d$ with trivial determinant. Suppose that there exists an $R$-linear surjection $\alpha:P\twoheadrightarrow I$, where $I\subset R$ is an ideal of height $d$. Then $P$ splits into a free summand of rank one if and only if $\mu(I)=d$.
\et
\proof First, we note that if $P$ splits into a free summand of rank one, then $\mu(I)=d$ follows using \cite[Lemma 1]{NMKL}. For a detailed proof, we refer to \cite[Corollary 4.4]{SMBB3}. Hence, we assume that $\mu(I)=d$. Let $I=\langle a_1,\ldots ,a_d \rangle$ be a set of generators for $I$. We fix an isomorphism $\chi: R\cong \wedge^d P$. Let the pair $\alpha,\chi$ induce $I=\langle b_1,\ldots,b_d \rangle+I^2$. We observe that in view of \cite[Corollary 3.4]{SMBB3}, to prove the theorem, it suffices to show that there exist $c_i\in I$ such that $I=\langle c_1,\ldots,c_d\rangle $ and $b_i-c_i\in I^2 $, $i=1,\ldots, d$. The remaining part of the proof is devoted to showing only this.

Considering that two sets of generators of $I/I^2$ may differ only by an invertible matrix in $R/I$, we find $\widetilde{\theta}\in \GL_{d}(R/I)$ connecting the two surjections $(R/I)^d\twoheadrightarrow I/I^2$ induced by $a_i$'s and $b_i$'s. We choose a lift $\theta \in \text{M}_{d\times d}(R)$ of $\widetilde{\theta}$ and an element $u\in R$ such that $\det (\theta)u-1\in I$. Now, we consider the unimodular row $w=(u, a_2, -a_1, a_3,\ldots,a_d)\in \Um_{d+1}(R)$. Since $\frac{1}{d!}\in R$ by Remark \ref{Rmk1} the unimodular row $w$ can be completed to the first row of a matrix in $\SL_{d+1}(R)$. Hence, using \cite[Proposition 7.4]{SBMKD2}, we can find a matrix $\tau\in\text{M}_{d\times d}(R)$ with $\det(\tau)=\det(\theta)$ modulo $I$, such that $(a_1,\ldots,a_d) \tau=( f_1,\ldots, f_d)$, where $I=\langle f_1,\ldots, f_d\rangle$.

Let ``bar'' denote going modulo $I$. It follows from the construction of $\tau$ that $\ol{\theta}^{-1}\ol{\tau}\in \SL_d(R/I)$. As $\dim(R/I)=0$ and $d\ge 3$, we have $\SL_d(R/I)=\E_{d}(R/I)$. Hence we can find a lift $\gamma\in \E_d(R)$ of $\ol{\theta}^{-1}\ol{\tau}$. Then the required $c_i's$ are defined as $(c_1,\ldots, c_{d}):= (f_1,\ldots, f_{d})\gamma$. This concludes the proof. \qed

\section{Cancellation of projective modules}\label{cansection}
This section is devoted to investigating the cancellation property of projective modules over a graded ring. We begin with a lemma, which is an analogy of \cite[Lemma 2]{P} in our setup.


\bl\label{cancri}
Let $R=\bigoplus_{i\ge 0} R_i$ and $M,M'$ be $R$-modules. Suppose that there exist $s,t\in R_0$ be co-maximal and isomorphisms $\sigma_1:M_s\iso M_s'$ and $\sigma_2:M_t\iso M_t'$ such that 
\begin{enumerate}
	\item $(\sigma_1)_t\equiv (\sigma_2)_s\mod(R_+)_{st}$;
	\item $M_{st}$ is a free module.
\end{enumerate}
Then there exists an isomorphism $\sigma:M\iso M'$ such that (i) $ \sigma_s\equiv \sigma_1\mod(R_+)_s$ and (ii) $\sigma_t\equiv \sigma_2\mod(R_+)_t$.
\el
\proof Let $\rank(M_{st})=n$ and let ``bar" denote going modulo $R_+$. Since $M_{st}$ is free there exists an isomorphism $\tau: M_{st}\iso R_{st}^{n}$. For an arbitrary isomorphism $\gamma:M_{st}\iso M_{st}$ we now consider the following commutative diagram
	$$\begin{tikzcd}
	M_{st}\arrow[r,"\gamma"]\arrow[d,"\tau"] &  M_{st} \arrow[d,"\tau"] \\
		R_{st}^{n}\arrow[r,"\widetilde{\gamma}"] & R_{st}^{n}
	\end{tikzcd}$$
where $\widetilde{\gamma}=\tau\gamma\tau^{-1}$. We will call $\widetilde{\gamma}$ is induced from $\gamma $ and $\tau$.

We take $\gamma=(\sigma_1)_{t}^{-1}\circ (\sigma_2)_s:M_{st}\iso M_{st}$. Then from (1) it follows that $\ol \gamma=\text{Id}$. We consider the isomorphism $\widetilde{\gamma}\in \GL_{n}(R_{st})$ induced from $\gamma$ and $\tau.$ Since $\ol \gamma=\text{Id}$, we have $\ol {\widetilde{\gamma}}=\text{Id}$. Applying Lemma \ref{isotopy} we can get a matrix, say $\widetilde{\theta(T)}\in \GL_{n}(R_{st}[T])$ such that $\widetilde{\theta(0)}=\text{Id}$ and $\widetilde{\theta(1)}=\widetilde{\gamma}$. Let us take $\theta(T)=(\tau\otimes R_{st}[T])^{-1}\widetilde{\theta(T)}(\tau\otimes R_{st}[T])$. Then we observe that $\theta(T)\in \Aut(M_{st}[T])$ such that $\theta(0)=\text{Id}$ and $\theta(1)=\gamma$. Now we define the isomorphism $\phi(T)=(\sigma_1\otimes R_{st}[T])\circ\theta(T):M_{st}[T]\iso M'_{st}[T]$. Then the proof follows from applying \cite[Lemma 1]{P}.\qed

Now we present the main theorem of the section.

\bt\label{can1}
Let $R$ and $d$ be as in Theorem \ref{umch}. Let $P$ be a projective $R$-module of rank $d$ such that $J(R_0,P)\not=0$. Then $P$ is cancellative. 

\et

\proof First we comment that since $J(R_0,P)\not=0$ there exists an $s\in R_0\setminus\{0\}$ such that $P_s$ is a free module. To see this let us choose a non-zero element $k\in R_0$ such that $P_k$ is an extended projective module from $R_0$. Consider the multiplicative set $T=R_0\setminus \{0\}$. Since $P_k$ is extended, the module $T^{-1}P_k$ is free. Now as $P$ is finitely generated we may choose a suitable multiple $s=kl$, for some $l\in T$ and ensure the existence of such an $s$.

 Let $(f,p)\in \Um(R\oplus P)$. Since $R$ has a non-trivial grading, the ideal $R_+$ is non-zero. In particular, the height of $R_+$ is $\ge1$. Hence, going modulo a non-zero element $g\in R_+$ and altering $(f,p)$ suitably via an element of $\text{Aut}(R\oplus P)$ we may assume that $f-1\in \langle g \rangle  $ and $p\in \langle g \rangle P$ [cf. Theorem \ref{umch}, first paragraph]. We take $P'=\frac{R\oplus P}{(f,p)R}$. Then to prove the theorem it is enough to show that there exists an isomorphism $\sigma:P\iso P'$. Moreover, we comment on an observation that finding a $\sigma:P\iso P'$ such that $\ol\sigma=\text{Id}$ is equivalent to find an $\alpha\in \Aut(R\oplus P)$ such that $\alpha(f,p)=(1,0)$ and $\ol \alpha=\text{Id}$.

 Let ``bar'' denote going modulo $R_+$ as well as $PR_+$. As $P_s$ is free, by Theorem \ref{umch} we can find $\alpha_1\in \Aut(R_s\oplus P_s)$ such that (1) $\alpha_1(f,p)=(1,0)$ and (2) $\ol \alpha_1=\text{Id}$. Then $\alpha_1$ will induce an isomorphism $\sigma_1:P_s\iso P'_t$ such that $\ol \sigma_1=\text{Id}$.
 
Let $S=\{1+sr:r\in R_0\}$. We denote $B=S^{-1}R$, $L=S^{-1}P$ and $L'=S^{-1}P'$. We note that $s\in \text{Jac}(S^{-1}R_0)$. Therefore, by Proposition \ref{dogr} we get $\dim(B_s)\le d-1$. Hence, using Lemma \ref{plgd} we can obtain a generalized dimension function $\delta:\Spec(B)\to \mathbb{N}$ such that $\delta(\p)\le d-1$ for all $ \p\in \Spec(B)$. Since $p\in \langle g \rangle P$, we note that $(f,p)\in \Um(R\oplus \langle g \rangle P)$. Moreover, the module $\langle g \rangle L$ is a projective $B$-module of rank $d$. Hence, applying Theorem \ref{eept} we can find $p'\in P$ such that $q:=p+gfp'\in \Um(\langle g \rangle L)$. Moreover, as $q\in \langle g \rangle L$ and $f-1\in \langle g \rangle $ one may obtain an $\alpha_2\in \text{Aut}(B\oplus L)$ such that $\alpha_2(f,p)=(1,0)$ and $\ol{\alpha}_2=\text{Id}$. Then $\alpha_2$ will induce an isomorphism $\sigma_2:L\iso L'$ such that $\ol{\sigma}_2=\text{Id}$. Since all modules are finitely generated (over a Noetherian ring) there exists an isomorphism $\sigma_2:P_t\iso P'_t$ such that $\ol{\sigma}_2=\text{Id}$, for some $t\in S$.

Now applying Lemma \ref{cancri} we get the required isomorphism $\sigma:P\iso P'$ such that $\ol \sigma=\text{Id}$. This completes the proof.\qed

\smallskip

\rmk Let $R$ and $P$ be as in Theorem \ref{can1}, and we consider $S=R_0\setminus \{0\}$. Then $J(R_0,P)\not=0 $ if and only if $S^{-1}P$ is free. If $S^{-1}R_0=\ol{\mathbb{Q}}$ and $\wedge^d P\cong R$, then it follows from \cite[Theorem 6.4.2]{AS} and \cite{S} that $J(R_0,P)\not=0$.
\smallskip

\rmk It would be interesting to know whether the hypothesis $J(R_0,P)\neq 0$ in Theorem \ref{can1} is necessary or not.




We now discuss an interesting consequence of Theorem \ref{can1}.

\bc\label{soprc}
	Let $R$ be an integral domain of dimension $d\ge 1$, and $A$ be a graded subring of $R[T]$ containing $R$ such that $\dim(A)=d+1$. Let $P$ be a projective $A$-module of rank $d+1$, so that the determinant of $P$ is extended from the base ring $R$. Then $P$ is cancellative.
\ec
\proof Let us take $\MT=R\setminus \{0\}$. Since $\dim(A)=d+1$, there exists an $a\in \MT$ such that $A_a$ contains a monic polynomial in $T$. Then $\MT^{-1}A\inj (\MT^{-1}R)[T]$ is an integral extension. This further implies that for any multiplicative set $S\subset R$ containing $a$, we must have $\dim(S^{-1}A)=\dim(S^{-1}R)+1$. Therefore, in view of Theorem \ref{can1} it is enough to show that $\hh(J(R,P))\ge 1$. To prove this we observe that, since $\dim(\MT^{-1}R)=0$, we have $\dim(\MT^{-1}A)= 1$. As the determinant of $P$ is extended from $R$, applying \cite{S} the module $\MT^{-1}P$ is a free $\MT^{-1}A$-module. Hence, there exists an element $s\in \MT$ such that $P_s$ is a free $A_s$-module. That is, the non-zero element $s\in  J(R,P)$.\qed

\section{Applications}\label{a1} This section is devoted to establishing some consequences of Theorem \ref{eg}.

\subsection{On a question of Nori: non-smooth graded case}
Let $R=\bigoplus_{i\ge 0} R_i$ be an affine graded domain of dimension $d\ge 3$ over a field $k$ such that $ \mathbb Q\subset k$ and $\dim(R_0)\ge 1$. In the next theorem, we show that Nori's question \cite{SM} on homotopy of sections of projective modules has an affirmative answer over $R$ even without the smoothness assumption.

\bt\label{NQ}
Let $R=\bigoplus_{i\ge 0} R_i$ be a graded domain of dimension $d\ge 3$ such that $\dim(R_0)\ge 1$. Moreover, we assume that $R$ is an affine algebra over an infinite field such that $\frac{1}{d!}\in R$. Let $I\subset R[T]$ an ideal such that $\mu(I/I^2T)=\hh(I)=d$. Then any set of generators of $I/I^2T$ lifts to a set of generators of $I$.
\et
\proof If $I$ contains a monic polynomial in $T$, then the result follows from \cite{SM}. Hence, without loss of generality, we may assume that $I$ does not contain a monic polynomial in $T$. Let $I=\langle f_1,\ldots,f_d\rangle +I^2T$. First, we comment that in \cite[Theorem 3.10]{MKD1} the hypothesis that the ring containing $\mathbb{Q}$ can be weakened by assuming the ring contains an infinite field such that $d!$ is invertible. We denote $R(T)=\MT^{-1}R[T]$, where $\MT$ be the ring consisting of all monic polynomials in $R[T]$. In view of \cite[Theorem 3.10]{MKD1} it is enough to prove that there exist $F_i\in IR(T)$, such that $IR(T)=\langle F_1,\ldots,F_d\rangle$ and $f_i-F_i\in I^2R(T)$. The proof is devoted to establishing only this.

Consider the multiplicative set $S=\{f\in R_0[T]:f\text{ is a monic polynomial}\}$ and let $B=S^{-1}R[T]$. We consider the grading $R[T]=\bigoplus_{i\ge 0}R_i[T]$. Let $\mm$ be a maximal ideal in $R_0[T]$. Then $\hh(\mm)=\dim(R_0)+1$. Therefore, from Suslin's monic polynomial theorem (see \cite[Chapter III, $\S 3$, 3.3, page no 108]{Lam}) we have $S\cap \mm\not=\emptyset$. Then by Proposition \ref{dogr} we have $\dim(B)\le d$. Since $S\subset R_0[T]$, the ring $B=\bigoplus_{i\ge 0} S^{-1}R_i[T] (=\bigoplus_{i\ge 0} B_i$ say) is also a graded domain of dimension $d$ such that $\dim(B_0)=\dim(S^{-1}R_0[T])=\dim(R_0)\ge 1$. As $I$ is not containing a monic polynomial we have $\hh(IB)\ge d$. Moreover, since $T$ is a unit in $B$, we have $IB=\langle f_1,\ldots,f_d\rangle B+I^2B$. Now, applying Theorem \ref{eg} we obtain $F_i\in IB$ such that $IB=\langle F_1,\ldots,F_d\rangle B$ and $f_i-F_i\in I^2B$. Since $B$ is a subring of $R(T)$ we get $IR(T)=\langle F_1,\ldots,F_d\rangle R(T)$ such that $f_i-F_i\in I^2R(T)$. This concludes the proof. \qed
\smallskip

\rmk\label{noh} It follows from \cite[Example 6.4]{BR} that in Theorem \ref{NQ} the hypothesis $\dim(R_0)>0$ is necessary.

\bc\label{egpa}
Let $R$ and $d$ be as in Theorem \ref{NQ}. Let $I\subset R[T]$ an ideal such that $\mu(I/I^2)=\hh(I)=d$. Then any set of generators of $I/I^2$ lifts to a set of generators of $I$.
\ec
\proof Let $I=\langle f_1,\ldots,f_d\rangle +I^2$ and let $I(0)=\{f(0):f\in R[T]\}$. Then as $R$ contains an infinite field without loss of generality we may assume that $I(0)\subset R$ is an ideal of height $\ge d$ (for details see the proof of \cite[Theorem 3.4]{BRS01}). If $\hh(I(0))>d$, then we can always lift any set of generators of $I(0)/I(0)^2$. Now if $\hh(I(0))=d$, then applying Theorem \ref{eg} there exist $a_i\in I$ such that $I(0)=\langle a_1,\ldots,a_d\rangle $, with $f_i(0)-a_i\in I(0)^2$, for $i=1,\ldots,d$. Hence, by \cite[Remark 3.9]{BR} there exist $g_i\in I$ such that $I=\langle g_1,\ldots,g_d\rangle +I^2T$ with $f_i-g_i\in I^2$ and $g_i(0)=a_i$, for $i=1,\ldots,d$. Now the result follows from Theorem \ref{NQ}.\qed

\bc\label{eueipe}
Let $R$ and $d$ be as in Theorem \ref{NQ}. Let $P$ be a projective $R[T]$-module with trivial determinant of rank $d$. Then $P$ has a unimodular element.
\ec
\proof Let us choose $\alpha\in P^*$. If $\alpha(P)=R[T]$, then there is nothing to prove. Hence without loss of generality we may assume that $\alpha(P)\subsetneq R[T]$. We apply Corollary \ref{ee} on the pair $(\alpha,1)$ to obtain an element $\beta\in P^*$ such that $\hh(I)\ge n$, where $I=(\alpha+\beta)(P)$. Now again, if $I=R[T]$, then the result is proved. Therefore, without loss of generality we may assume that $I$ is a proper ideal in $R[T]$. Then again by Corollary \ref{ee} we get that $\hh(I)=d$. Now the result follows from Corollary \ref{egpa} and subtraction principle as stated in \cite[Corollary 4.13]{MKD1} (taking $Q=(R[T])^{d-1}$, $I_1=R[T]$ and $I_2=I$).\qed
\smallskip

\rmk One can remove the restriction on the base field in Corollary \ref{eueipe} in the following way: let $P$ be a projective $R[T]$-module of rank $d$ with trivial determinant. Recall that the ring $R(T)$ is obtained by localizing $R[T]$ with respect to the multiplicative set consisting of all monic polynomials in $R[T]$. Then, in view of \cite[Theorem 5.2 and Remark 5.3]{SMBHLRR}, it is enough to show that the modules $P/TP$ and $P\otimes R(T)$ have unimodular elements. Let $S$ be the multiplicative set consisting of all monic polynomials in $R_0[T]$. Then, it follows from Corollary \ref{eue} that $P/TP$ and $S^{-1}P$ (and hence $P\otimes R(T)$) have unimodular elements.
\smallskip


	\subsection{Generating ideals up to projective equivalence}
	Recall that, two ideals $I \text{ and } J$ in a ring $ A$ are said to be \textit{projectively equivalent} if some power of $I$ and some power (usually different) of $J$ have the same integral closure. The following theorem is an improvement of \cite{katz} in our setup.

\bt\label{pe}
Let $A$ and $n$ be as in Theorem \ref{eg}. Let $I\subset A$ be an ideal of height $\ge 2$. Then there exists an ideal $J\subset A$ projectively equivalent to $I$ satisfying $\mu(J)\le n$.
\et
\proof First we observe that combining the results \cite[Proposition 2.2]{DRA} and Theorem \ref{eg} one can prove the following: let $K\subset A$ be an ideal such that (i) $\mu(K/K^2)\le n$ and (ii) $\hh(K)\ge 2$. Then $\mu(K)\le n$. Applying \cite{katz} we obtain an ideal $J\subset A$ such that (1) $I$ and $J$ are projectively equivalent, (2) $\hh(J)\ge 2$ and (3) $\mu(J/J^2)\le n$. Now it follows from the previously mentioned observation that $\mu(J)\le n$.\qed

\section*{Acknowledgment}  We thank Mrinal Kanti Das for suggesting Theorem \ref{pe}. The author is grateful to the referee(s) for their careful reading and valuable suggestions. Without their detailed comments, the exposition would have lacked clarity at various places.

         \bibliographystyle{abbrvurl}

		\end{document}